\documentstyle{amsppt}

\NoBlackBoxes


\def\cal{\Cal}


\def\zx#1{{\tenbf #1}}
\def\xbb{Let $X$ be a Banach space.  }


\def\rptens#1#2{${\smash{\mathop{\widehat{\otimes}}}}^{#1} #2$}
\def\ptens#1#2{${\smash{\mathop{\widehat{\otimes}}}}^{#1}_{s} #2$}

\def\mpotimes{\smash{\mathop{\widehat{\otimes}}}}

\def\mrptens#1#2{{\smash{\mathop{\widehat{\otimes}}}}^{#1} #2}
\def\mptens#1#2{{\smash{\mathop{\widehat{\otimes}}}}^{#1}_{s} #2}

\def\mtens#1#2{\smash{\mathop{\widehat{\otimes}}\limits_{#1} } #2}

\def\infseq#1#2{${\{{#1}_{#2}\}}_{{#2}=1}^\infty$}
\def\minfseq#1#2{{\{{#1}_{#2}\}}_{{#2}=1}^\infty}

\def\psh{${\cal P}-$Schur}
\def\pnsh{${\cal P}_N-$Schur}

\def\pnwe{${\cal P}_N-$weak}

\def\pnwey{${\cal P}_N-$weakly}
\def\diag#1#2{${\Delta}_{#1} (#2)$}

\def\mnm#1{{\|#1\|}}

\def\mdiag#1#2{{\Delta}_{#1} (#2)}

\def\pol#1#2{${\cal P}_{#1}{(#2)}$}

\def\mpol#1#2{{\cal P}_{#1}{(#2)}}
\def\pndp{${\cal P}_N$ Dunford--Pettis}
\def\pndpp{${\cal P}_N$ Dunford--Pettis property}

\leftheadtext{JEFF FARMER AND WILLIAM B. JOHNSON}
\rightheadtext{POLYNOMIAL SCHUR AND POLYNOMIAL DUNFORD--PETTIS
PROPERTIES}

\topmatter
\title Polynomial Schur and\\
Polynomial Dunford-Pettis Properties\endtitle
\author Jeff Farmer and  William B. Johnson\endauthor

\address Department of Mathematics, Texas A\& M University,
College Station, Texas 77843\endaddress

\curraddr Department of Mathematics, University of Missouri,
Columbia, Missouri  ZZZZIPPPP\endcurraddr

\email\endemail

\address Department of Mathematics, Texas A\& M University,
College Station, Texas 77843
\endaddress

\email wbj7835@tamvenus, wbj7835@venus.tamu.edu\endemail


\abstract
A Banach space is {\it polynomially Schur}  

if sequential convergence against analytic polynomials implies norm
convergence.  Carne, Cole and Gamelin show that a space has this  
property
and the Dunford-Pettis property if and only if it is Schur.  Herein  
is defined a
reasonable generalization of the Dunford--Pettis
property using polynomials of a fixed
homogeneity.  It is shown, for example,  that  

a Banach space will has the $P_N$ Dunford--Pettis property
if and only if every weakly compact $N-$homogeneous
polynomial (in the sense of Ryan) on the space is completely  
continuous.
A certain geometric condition, involving estimates on spreading  
models and 

implied by nontrivial type, 

is shown to be sufficient to imply that a space is
polynomially Schur.

\endabstract

\thanks This paper forms a portion of the first author's doctoral
dissertation written under the supervision of the second author.
\endgraf The first author was supported in part by NSF Grant  
\#DMS-9021369.
\endgraf The second author was supported in part by NSF Grant  
\#DMS-9003550.
\endgraf
This paper is in final form and no version of it will be submitted  
for
publication elsewhere.\endthanks
\endtopmatter

\document
\head 1.  Introduction
\endhead
The relationship between holomorphic functions defined on an infinite
dimensional
Banach space and (geometric or topological) properties
of the space has been of recent interest (see, for example,  
\cite{AAD},
\cite{ACG},  \cite{CCG},  \cite{CGJ},  \cite{F},  \cite{R 1}).  As in  
the
one--dimensional case, holomorphic functions are defined in terms of  
Taylor
series, which in the
infinite--dimensional case have terms consisting of homogeneous  
analytic
polynomials.
Just as in the case of linear functionals (1--homogeneous  
polynomials),
one can
consider properties of the topologies
induced by the polynomials on the space.
In this paper we consider the properties which are
analagous to the Schur property
and the Dunford--Pettis
property; i.e., those obtained by replacing weak
sequential convergence
with sequential convergence against an
arbitrary $N$--homogeneous analytic polynomial.  We
relate these properties to one another and to the geometric
property of type
and the existence of certain 

spreading models.
\par
$X$ will be a complex infinite--dimensional Banach
space.  An {\it $N-$homogeneous
analytic polynomial} on $X$ is the
restriction to the diagonal of an $N-$linear form on the $N-$fold
Cartesian product of $X$ with itself, or equivalently, a linear  
functional
on the
$N-$fold projective tensor product of $X$ with itself.  Indeed, given  
an
$N-$homogeneous analytic function $P$ on $X$, one obtains an  
$N$--linear
form on $X$ by taking the $N$th derivative and dividing by $N!$; the  
form is
related to the polynomial by the polarization formula:
$$A_P(x_1,\dots ,x_n)=\text{Avg}\{{\epsilon}_i=\pm 1\}(\Pi  
{\epsilon}_i)
P\bigl(\sum_{i=1}^n{\epsilon}_i x_i\bigr).$$

The form $A_P$ is clearly symmetric (invariant under permutations of  
the
coordinates).  Likewise any bounded symmetric $N-$linear form will  
give
rise to an $N-$homogeneous analytic polynomial. Such a form can be  
linearized
by taking the projective tensor product of $X$ with itself $N$
times and extending 

the form to a linear functional on this tensor product.
The subspace of symmetric
linear functionals is the dual of the symmetric
$N-$fold projective tensor product,
which is a complemented subspace of the $N-$fold projective tensor  
product.
The projection is given by extending the following map linearly:

$$(x_1\otimes x_2\otimes\dots\otimes x_n)\to {1\over  
{n!}}\sum_{\pi\in S_n}
(x_{\pi_1}\otimes x_{\pi_2}\otimes\dots\otimes x_{\pi_n}).$$
 We denote the symmetric projective tensor product by  \ptens N X .   
The
$N$--linear
form $A_P$ associated with $P$ can now be considered a linear  
functional on
\ptens N X . The supremum norm of the polynomial is related to that  
of the
linear functional as follows:
$$||P||\le||A_P||\le {{N^N}\over{N!}} ||P||$$
If we call the space of polynomials
${\cal P}_N$ the above simply says that $ {\cal P}_N$ is isomorphic  
to
${(\mptens N X)}^*$.
Since for our purposes the index $N$ will be fixed,
we will suppress reference to
this isomorphism and use the same label for a polynomial and its  
associated
symmetric linear functional.  More details about the above  
relationships may be
obtained from \cite{M} or \cite{R 1}.
 We will study the topologies generated by these
polynomials, especially with respect to sequential convergence.
\par We define the \pnwe\ topology on X to be the topology generated
by the all of the homogeneous analytic
polynomials of degree less than or equal to
$N$; that is, a net $\{x_\alpha\}$
converges to $x$ in the  $ {\cal P}_N-$weak topology if, for every  
$M\le N$,
for every $M-$homogeneous
analytic polynomial $P$, $P(x_\alpha )\to P(x)$.  Note that for $N=1$
this is the usual weak topology, and that for $M>N$ the $ {\cal  
P}_M-$weak
topology is finer than the \pnwe\ one.  We call the weak polynomial
topology the topology
which is generated by the union of $ {\cal P}_N$ for all $N\in \zx  
{{Z}_+}$.
 In analogy to the Schur property, we say a space is
{\it \pnsh}\ if whenever $P(x_n)\rightarrow 0$ for all $P \in {\cal  
P}_N$
then $x_n$ is norm null.  If  $P(x_n)\rightarrow 0$ for all $P \in  
{\cal P}_N$
for all $N$ implies that $x_n$ is norm null, then we say $X$ is {\it  
\psh}.
It is evident (multiply linear functionals) that every \pnsh\ space
is \psh\ and that every Schur space is \pnsh\ for every $N$ (and  
\psh).
\par These topologies were introduced in  \cite{R 1}, and the weak  
polynomial
topology also appeared in  \cite{CCG} which considered relations  
between the
Dunford-Pettis property, the Schur property and the \psh\ property  
(in the
terminology of  \cite{CCG}, ``$X$ is \psh\ "$=$ ``$X$ is a $\Lambda  
-$space'').
\smallskip\par
Let $\theta :X\to \mptens N X $
by $\theta (x) = x\otimes \cdots\otimes x$ (N times)
and define $\theta (X) =$\diag  N X.
\smallskip
 This set is a non-convex,
norm-closed subset of $\mptens N X$ with
the property that $\lambda x\in \mdiag N
X$ whenever  $x\in \mdiag N X$.
\smallskip
Now $\theta$ is a continuous $N-$homogeneous
function, which is the (nonlinear) preadjoint of the isomorphism  
between
$ {\cal P}_N$ and ${(\mptens N X)}^*$.
Notice that
$\theta$ is an $N-$to$-$one map; we have
$$ {\theta}^{-1}(x\otimes\cdots\otimes x)=\{e^{{2\pi  ni}\over N}  
x|n\in \zx Z
\}$$
(Use separating functionals to the $N$th power to show equality.)
We reserve the symbol ${\theta}_N$ for this function.
\par
If $P(x_\alpha)\rightarrow P(x)$ for all $P\in \mpol N X$ then the  
net
need not converge against polynomials in ${\cal P}_M$ for all $M<N$,  
but since
$P(x)=P(y)$ for all $P \in {\cal P}_N$
implies that $x\over y$ is a complex $N$th
root of unity, if also $x_\alpha\rightarrow x$ weakly, then  
$x_\alpha$ 

${\cal P}_N-$weakly converges to $x$.  Thus in practice it is easy to  
pass
from convergence against all $N-$homogeneous polynomials to $ {\cal  
P}_N-$weak
convergence.  \smallskip\smallskip
 Although for any one polynomial $P$, $P(x- x_{\alpha})\to 0$ and
$P(x_{\alpha})\to P(x)$ are not in general equivalent, the following
known fact  is useful.
\proclaim  {Lemma 1.1}
 $x_{\alpha}\to x$ in the \pnwe\ topology, if and only if $x-  
x_{\alpha}\to 0$
in the \pnwe\ topology.
\endproclaim
\demo{Proof (sketch)}
Let $P$ be an $N$-homogeneous polynomial and let $x_{\alpha}\to x$ in  
the \pnwe
\ topology.  Then, letting $A_P$ be the $N$--linear form associated  
with $P$ we
have

$$P(x- x_{\alpha})=\sum_{i=1}^N {(-1)}^i {N\choose i}
{A_P}(x,x,\dots,x,x_{\alpha},\dots,x_{\alpha})$$
where in each term $x$ appears $i$ times and $x_{\alpha}$ appears  
$N-i$ times.
Since convergence in the \pnwe\ topology implies convergence against  
any
polynomial of lesser homogeneity, we consider each term as
an \hbox{$(N-i)-$}homogeneous polynomial ($x$ being fixed), to see  
that
the sum indeed converges to zero.  The converse is obtained, using  
the same
expansion,  by induction on $N$.
\enddemo

\head 2.  Polynomial Dunford--Pettis Spaces
\endhead
\par One result of  \cite{CCG} is
that a Banach space is Schur if and only if it is
polynomially Schur and has the Dunford--Pettis property.  We can  
obtain an
analagous result for polynomials of fixed homogeneity by
defining an appropriately
analagous Dunford-Pettis property.
Our first task is to adapt Lemma 7.3 of  \cite{CCG} for our purposes.

\proclaim{Propopsition 2.1}  The following are 

equivalent  for any Banach space $X$, and any
fixed positive integer $N$.
\roster
\item"(i)"Any polynomial on $X$ is \pnwey\ sequentially continuous.

\item"(ii)"If \infseq x n  is a \pnwey\ null sequence in $X$ (i.e. if
$\{x_k\otimes\cdots\otimes x_k\}$ \hbox{(N times)} is weakly null),  
then
$\{x_k\otimes\cdots\otimes x_k\}$ (m times) is weakly null in
$\mrptens m X$ for $m>N$.

 \item"(iii)"For $m>N$ the function $\theta =\theta (N,m)$
which takes
$\theta :\mdiag N X \to \mdiag m X$ by
$$
x\otimes\cdots\otimes x \enspace ({\tenrm N \enspace times})\mapsto
x\otimes\cdots\otimes x\enspace ({\tenrm m\enspace times})
$$
is weak to weak sequentially continuous.
\endroster
\endproclaim
\demo{Proof}
The proof of these equivalences is an exercise, and can be
adapted easily from the
proof of Lemma 7.3 in  \cite{CCG}.
 \enddemo
\par  For $n=1$, the equivalent properties of 2.1 were shown to
be implied by
the Dunford--Pettis property.   We will now define a polynomial
Dunford--Pettis property
which will imply the conditions of 2.1 for each positive integer.
\par  We say that a space $X$ has the \pndpp\ provided that it  
satisfies any
of the equivalent conditions of Proposition 2.2.

\proclaim {Proposition 2.2}
 \xbb For fixed $N$, the following are equivalent :
\roster
 \item"(i)"Whenever $\minfseq P n$ is a weakly null sequence of  
$N-$homogeneous
 polynomials
(or equivalently, symmetric bounded $N-$linear forms on $X$) and
\infseq x n converges \pnwey\ to $x$ in $X$, then $P_n(x_n)\to 0$.

 \item"(ii)" Every weakly compact operator on $\mptens N X$
is completely
continuous when restricted to \diag N X .

 \item"(iii)"If K is a weakly compact set in any Banach space Y, and
$J$
 is \pnwey\ compact in
$X$, then  ${\theta}_N(J)\otimes K$ is a weakly compact set in
$(\mptens N X )\mrptens {} Y$
\endroster
\endproclaim
\noindent where by  ${\theta}_N(J)\otimes K$ we mean simply the set
${\theta}_N(J)\times K$, that is, the set of all
 ${\theta}_N(j)\otimes k$ with $j\in J$ and $k\in K$.

\par  In the case $N=1$ these conditions reduce to known equivalent  
statements
of the classical Dunford--Pettis property; this proposition justifies 

the definition of the \pndpp\.
Before giving the proof, we make the following remark.
\par
 R. Ryan in  \cite{R 2} considers $N-$homogeneous polynomials from  
$X\to Y$;
as in the scalar case, we can equivalently consider linear operators  
from
\ptens N X to $Y$; such a polynomial is weakly
compact if it maps bounded sets to weakly compact ones (i.e. if the  
associated
linear operator is weakly compact).
Ryan investigates some conditions which are
equivalent to weak compactness of such polynomials.
Using this definition it is easy to see that (ii) above is equivalent  
to
\proclaim{}
\item{${{\tenrm(ii)}}^\prime$} Every weakly compact $N$--homogeneous
polynomial from $X$ to any Banach space $Y$ is completely continuous  
(on $X$).
\endproclaim

\demo{Proof}
\roster
\item"(i)$\Rightarrow$(iii)"
We want to show that  ${\theta}_N(J)\times K$
is weakly compact in $\mptens N X\otimes Y$.
Take a sequence ${\theta}_N(x_n)\otimes k_n$ in  $  
{\theta}_N(J)\otimes K$; by
hypothesis assume we have passed to
 a subsequence such that
 ${\theta}_N(x_n)$ and $k_n$ are weakly convergent to
 ${\theta}_N(x)$ and $k$ in \ptens N X and  $Y$, respectively.
That is,  ${\theta}_N(x_n) - {\theta}_N(x)$ and $k_n -k$ are weakly  
null.
If
$$\phi\in {(\mptens N X \otimes Y)}^*\equiv {\cal B}
\bigl(Y,{(\mptens N X )}^*\bigr)$$
then $\phi (k_n -k)$ is weakly null by continuity and (i) applies.
We then have
$$\Bigl<\phi (k_n -k),{\theta}_N(x_n) - {\theta}_N(x)\Bigr> =$$
$$\qquad\bigl<\phi (k_n),{\theta}_N(x_n) \bigr> -
\bigl<\phi (k),{\theta}_N(x_n)\bigr> -\bigl<\phi  
(k_n),{\theta}_N(x)\bigr> +
\bigl<\phi (k),{\theta}_N(x)\bigr> $$
 Taking limits as $n\to\infty$ we see
$$\quad 0={
\lim\limits_{n\to\infty}\bigl<\phi (k_n),{\theta}_N(x_n) \bigr> -2
\bigl<\phi (k),{\theta}_N(x)\bigr>  +
\bigl<\phi (k),{\theta}_N(x)\bigr>}$$
$$=
{\lim\limits_{n\to\infty}\bigl<\phi (k_n),{\theta}_N(x_n) \bigr> -
\bigl<\phi (k),{\theta}_N(x)\bigr> }
$$
This says exactly that
 ${\theta}_N(x_n)\otimes k_n$ is weakly convergent to  
${\theta}_N(x)\otimes k$.
\smallskip

\item"(iii)$\Rightarrow$(ii)"
Let $T:\mptens N X \to Y$ be weakly compact
 and ${\theta}_N(x_n)\to {\theta}_N(x)$ weakly in the symmetric  
tensor product.
Choose ${\phi}_n$ to be norming functionals for
$T({\theta}_N(x_n)-{\theta}_N(x))$
in the sphere of
$Y^*$.  Since $T^*$ is also weakly compact, assume by passing to a  
subsequence
that $T^*({\phi}_n)$ converges weakly,
say to $\psi$.
By passing to a subsequence we can also assume that
\hbox{$\{{T^*}({\phi}_n)-\psi\}_{n=1}^{\infty}$} is basic (or norm  
null, in
which case the argument is simpler).
Apply (iii) to the sets

$$K={\bigl\{T^*({\phi}_n)-\psi\bigr\}}_{n=1}^\infty\cup \{0\}\quad
{\tenrm and}\quad J=\minfseq x n \cup\{x\}$$
to see that
 ${\theta}_N(J)\otimes K$   is weakly compact in
$(\mptens N X )\mtens {} {(\mptens N X)}^*$.
The sequence
 ${\theta}_N(x_n) \otimes\{T^*({\phi}_n)-\psi\} $ thus has
a convergent subsequence
(we pass to that).
 First we claim that this subsequence must go weakly
to zero; indeed, it goes to zero in a weaker Hausdorff topology,  
namely that
generated by considering the weak topology on the second co-ordinate.
Since it is clear that ${\theta}_N(x) \otimes\{T^*({\phi}_n)-\psi\}$
is weakly null, we can conclude that
$$\quad w-\lim \Bigl([{\theta}_N(x_n)-{\theta}_N(x)]
\otimes\{T^*({\phi}_n)-\psi\}\Bigr)=0$$
Now consider the functional associated with the identity operator;  
call it
$\Gamma$. We have
$$\eqalign{\quad 0&={\lim\limits_{n\to\infty} \Gamma
\bigl(\{T^*({\phi}_n)-\psi\}\otimes  
[{\theta}_N(x_n)-{\theta}_N(x)]\bigr)}\cr
\quad&={\lim\limits_{n\to\infty}
\bigl<\{T^*({\phi}_n)-\psi\}, {\theta}_N(x_n)-{\theta}_N(x)\bigr>}\cr
\quad&={\lim\limits_{n\to\infty}
\biggl(\bigl<T^*({\phi}_n),  
{\theta}_N(x_n)-{\theta}_N(x)\bigr>-\bigl<\psi,
{\theta}_N(x_n)-{\theta}_N(x)\bigr>}\biggr)\cr
\quad
&={\lim\limits_{n\to\infty}
\bigl<{\phi}_n, T({\theta}_N(x_n))-T({\theta}_N(x))\bigr>}\cr
&={\lim\limits_{n\to\infty}
\mnm{T({\theta}_N(x_n))-T({\theta}_N(x))}}}$$
This gives (ii).

\item"(ii)$\Rightarrow$(i)"
Let \infseq P n be the weakly null sequence  and \infseq x n go  
\pnwey\ to $x$
 and define a map
$T$ from \diag N X to $c_0$ by
$$T({\theta}_N(z))={\bigl(P_n(z)\bigr)}_n\quad\forall z\in X$$
The map $T$ extends linearly (via the polarization formula) to all
of \ptens N X. Since $T^*(e_n)=P_n$ goes weakly to 0 we see that the  
map $T$
is weakly
compact.  Applying (ii) we get $T({\theta}_N(x_n))$ going in the norm  
on
$c_0$ to  $T({\theta}_N(x))$. But since the norm on $c_0$ is the sup  
norm,
this gives (i) and completes the proof.
\endroster
\enddemo
We note that the condition (ii) is sharply stated with the following  
example.
\example{Example}  We will see momentarily that $l_2$ is ${\cal  
P}_2$--Schur
and therefore
is ${\cal P}_2$--Dunford Pettis.  Consider the operator
$$Id{\otimes} Q_1:l_2{\mpotimes} l_2\to c_0$$
where $Q_1$ is the projection onto the first basis vector. Consider  
the
symmetrized version, that is, restrict the operator to the symmetric
tensor product, which is a complemented subspace.  This operator is  
weakly
compact and therefore completely continuous on ${\theta}_2 (l_2)$ by
Proposition 2.2 but is clearly not
completely continuous on the entire symmetric
tensor product; consider the image of $e_n\otimes e_1+ e_1\otimes  
e_n$, for
example, which is weakly null but whose image is the unit vector  
basis of $c_0$.
\endexample
\smallskip
\proclaim{Proposition 2.3}
 If a Banach space is \pndp\ then it satisfies the equivalent
conditions in Proposition 2.1.
\endproclaim
\demo{Proof}
We prove 2.2(i) implies 2.1(ii).\smallskip
Let \infseq x n be a \pnwey\ null sequence in $X$, i.e.
${{\theta}_N(x_n)}_{n=1}^\infty$ is weakly null in
\rptens N X and
${{\theta}_M(x_n)}_{n=1}^\infty$ is also weakly null in \rptens M X  
whenever
\hbox{$1\le M<N$}.
Let $m=N+1$, $\phi\in {(\mrptens {N+1} X )}^*$ and consider
$\phi$ as a linear operator from $X$ to ${(\mrptens N X)}^*$.
Since \infseq x
n is weakly null in $X$, so is its image in ${(\mrptens N X)}^*$  
under $\phi$.
Thus
$${\bigl<\phi ,{\theta}_{N+1}(x_n)\bigr>}=
{\bigl<\phi (x_n) ,{\theta}_{N}(x_n)\bigr>}\to 0$$
by the first formulation of the {\pndpp} (notice that for this  
application it
matters not whether $\phi$ is symmetric).  This proves the  
proposition for
$m=N+1$ and by induction (in an obvious way) for $m=qN+1$ for  
$q=1,2,\dots$.
But we can also write an analogous proof for $m=N+k$ for $2\le k<N$  
and extend
it by induction as well.\enddemo

\proclaim{Proposition 2.4}
\xbb  For fixed $N$, the following are equivalent :
\roster
\item"(i)" $X$ is {\pnsh}.

\item"(ii)" $X$ has the {\pndpp}\ and is {\psh}.

\item"(iii)" $X$ satisfies (i)--(iii) of proposition 1.1 and is  
{\psh}.
\endroster
\endproclaim
\demo{Proof}
(i)$\Rightarrow$(ii)  requires only Lemma 1.1 and
(ii)$\Rightarrow$(iii) is Proposition 2.3, so
(iii)$\Rightarrow$(i) remains.
Let ${\theta}_N (x_n)$ be weakly null.  Then  ${\theta}_M (x_n)$
is weakly null for all $M$ by 2.1.  But since $X$ is {\pnsh},
$x_n$ must go to 0 in norm. \enddemo

It is of interest to note that if we are not in the context of the  
\psh
\ property, the conditions of 2.1
are weaker than the \pndpp; $T^*$, the original Tsirelson
space (or, in fact any space having the approximation property
 with \pol N X reflexive for all $N$;
see [F] for further discussion of such spaces) will satisfy 2.1 for  
all $N$
but fail to be \pndp.
\example{Examples}
It is clear from the classical work of Pitt  \cite{P} that $l_p$  
spaces for
($1\le p<\infty$) are \pnsh\ for $N\ge p$, and it is proved in  
\cite{CCG}
that $L_p$ spaces  ($2\le p<\infty$) are \psh\ (in fact \pnsh\ for  
$N\ge p$);
we can thus conclude that they are \pndp.  The space $c_0$ is  
Dunford--Pettis
and therefore {\pndp}\ for every $N$.  This implies  (for example)  
that
$l_3\oplus c_0$ is ${\cal P}_3$--Dunford--Pettis but not \pnsh\ for  
any $N$.
In the next section we discuss further exactly which spaces may be  
\pnsh.
\endexample

\head 3.  Spaces with Type are Polynomially Schur
\endhead

\bigskip
\par
In this section we will give some sufficient criteria for spaces to  
be
Polynomially Schur. We will show, for example, that any space
having non-trivial type is \psh\ and indeed is \pnsh\ for some $N$.  

(Jaramillo and Prieto  \cite{JP} have independently shown that every
superreflexive space is polynomially Schur).  In particular, $L_p$  
spaces 

are Polynomially Schur for all $1 < p < \infty$.
\par
It is convenient to use the concept of a spreading model, the
construction of which is due to Brunel and Sucheston  \cite{BS 1}.  
\par
Finite versions of Ramsey's Theorem allow that given any property of
$n$-tuples of elements from a sequence, one can pass to a subsequence  
with the
property that all $n$-tuples formed from the subsequence share the
property or
else all fail it.  By using the size of the norm of a sum of $n$  
elements as
the property one can, by repeatedly
applyling the theorem, approximately stabilize
the norm (to within any desired ${\epsilon}_n$) of any finite  
combination as
long as many of the beginning terms are thrown away.  More preciesely  
we have
the following fact (see [B] or [BS 1]):

\proclaim{Proposition 3.1} Let $(f_n)$ be a bounded sequence with no
norm-Cauchy subsequence in
a Banach space $X$.  Then there exists a subsequence $(e_n)$ of   
$(x_n)$ and a
norm $L$ on the vector space $S$ of finite sequences of scalars such  
that
$$\forall \epsilon>0\quad\forall a\in S \quad\exists k\in N\enspace
\text{s.t.}\enspace\forall k<k_1<k_2<\dots <k_M$$
we have
$$
\biggl|\Bigl\Vert\sum a_i e_{k_i}\Bigr\Vert-L(a)\biggr|<\epsilon
$$
\endproclaim
The completion of $[e_i]$ (call it $F$)
under the norm $L$ is called a spreading model for the sequence  
$(e_n)$.  The
reason for the terminology is that the sequence $(e_n)$
is invariant under spreading with respect to the norm $F$, that is to  
say,
for every finite sequence of scalars $(a_i)$ and every subsequence  
$\sigma$
of the natural numbers
$$
 {\biggl\Vert \sum_{i=1}^{M} a_i
e_i\biggr\Vert}_F= {\biggl\Vert \sum_{i=1}^{M} a_i
e_{\sigma(i)}\biggr\Vert}_F
$$
Thus any norm estimate satisfied by sequences in the spreading model  
will be
approximately satisfied for sequences of finite length to any desired  
degree
provided we go out far enough in the sequence $(x_n)$.
If the original sequence
was weakly null then the resulting sequence will be
unconditional; that is to say, we have the following
(Lemma 2 of  \cite{B}, or see  \cite{BS 2}).
\proclaim {Proposition 3.2}
 If $(x_n)$ is weakly null, then the sequence $(e_n)$ is  
unconditional
in $F$ with
unconditional constant at most 2.\endproclaim
\medskip
Now we are ready to state the criterion.
\smallskip
\proclaim{Theorem 3.3}
Suppose a Banach space $X$ has the property that for every normalized
weakly null sequence $\{y_n\}$ in $X$ there exists a subsequence and  
a sequence
 $\{f_n\}$ in $X^*$ biorthogonal to it which has an (unconditional)
spreading model with an
upper $p$-estimate for some $p>1$. Then $X$ is \psh.   If the same  
$p$
works for every such sequence and $N>p^\prime$,
then the space $X$ is \pnsh.
\endproclaim
\par
We know of no space which is \psh\ but which
fails the above property.  The space
$
{(l_3{\oplus} l_4{\oplus} l_5{\oplus}\cdots)}_2
$
is easily seen to be \psh\ although it fails cotype (and hence type  
and
superreflexivity), yet is reflexive; it does satisfy the hypothesis  
of 3.3.
A Schur space satisfies the hypothesis vacuously.
\demo{Proof}
We will pass to subsequences and relabel without mercy.  Start with  
any bounded
sequence in $X$ which is not norm null; we need to find a polynomial  
which is
bounded away from zero on a subsequence.  By Rosenthal's theorem  
\cite{D\rm,
Chapter XI} there is either a weakly Cauchy subsequence or a  
subsequence
equivalent to the unit vector basis of $l_1$.  Since the unit vector  
basis of
$l_1$ is not weakly null, we are done in this case.  For the same  
reason
(linear
functionals are polynomials) we are finished if our weakly Cauchy  
sequence is
not weakly null.  So we have reduced to the case of a bounded weakly  
null
sequence
which is not norm null and can apply the 

hypothesis to that sequence 

(it is purely formal that the ``normalized" 

condition can be replaced by ``bounded").
\par
Let $ \{y_n,f_n\}$ be the biorthogonal system obtained, and assume  
with
no loss of generality that the
$f_n$ are normalized.  By the 

definition of a spreading model we
know that
for any $c>0$ we can find a constant $C$ so that for every $M$ we  
have
$$
{\biggl\Vert \sum_{i=c\log M}^{M}
f_i\biggr\Vert}\le C {\biggl(\sum_{i=c\log M}^{M}
{\Vert f_i\Vert}^p\biggr)}^{1\over p}.
$$\smallskip\noindent
This is because once we have a spreading model we may always improve  
the
stability estimates by passing to a subsequence.
 Choosing $c$ small enough so that $m=c\log M\le M^{1\over p}$ for  
all
$M\in{\cal N}$, and letting $(a_i)$ be scalars of modulus less than  
or
equal to 1, we obtain

$$
{\biggl\Vert \sum_{i=1}^{M}
a_if_i\biggr\Vert}\le \sum_{i=1}^{m-1}||f_i|| +
{\biggl\Vert \sum_{i=m}^{M}
a_if_i\biggr\Vert}\le  m + C{\biggl(\sum_{i=m}^M
{ |a_i|}^p\biggr)}^{1\over p}\le C_1M^{1\over p}.
$$ \smallskip\noindent
Now by general Banach lattice techniques (see Lemma 3.4 below), we
know that for any $r<p$ (we choose $r$ so
that $p^\prime<r^\prime<N$) we can get a different $C$ so that
for any sequence $(a_i)$ we will have
$$
{\biggl\Vert \sum_{i=1}^{M} a_i
f_i\biggr\Vert}\le C{\biggl(\sum_{i=1}^M {|a_i|}^r\biggr)}^{1\over  
r}.
$$\smallskip\noindent
Now, given $y\in B$ choose  $(a_i)\in l_r$ of $l_r$ norm one
to norm the sequence ${(f_i(y))}_{i=1}^M$ in
$l_{r^\prime}$.  We then have that
$$
{\Biggl|\biggl(\sum_{i=1}^M {f_i(y)}^N\biggr)\Biggr|}^{1\over N}
\le {{\biggl(\sum_{i=1}^M {|a_i|}^r\biggr)}^{1\over r}
{\biggl(\sum_{i=1}^M {|f_i(y)|}^{r^\prime}\biggr)}^{1\over  
{r^\prime}}}.
={\sum_{i=1}^M a_if_i(y)}
$$

$$\le {\biggl\Vert \sum_{i=1}^{M} a_i
f_i\biggr\Vert}\le C{\biggl(\sum_{i=1}^M {|a_i|}^r\biggr)}^{1\over r}
=C
$$ which says that

$$
P(y)=\sum_{i=1}^\infty {f_i(y)}^N \le C^N
\quad\text{hence}\quad ||P||\le C^N .$$
But we know that
$P(y_n)={f_n(y_n)}^N$ is bounded away from zero because the
sequence $(y_n,f_n)$ was biorthogonal.
Therefore no sequence bounded away from 0 in $X$ can be polynomially
null, and
so $X$ is \psh.  If there is a uniform value for $p$ then $X$ is
\pnsh\ for $N>p^\prime$.
\enddemo
\par
It remains for us to prove the following standard fact.

\proclaim {Lemma 3.4}
Suppose that $\{x_i\}$ is a normalized sequence in a Banach space  
satisfying
$$
{\biggl\Vert \sum_{i\in B}
a_ix_i\biggr\Vert}
\le C|B|^{1\over p}\max_{i\in B}|a_i|
$$
for all scalars $(a_i)$, and all finite subsets $B$ of the natural  
numbers.
Then for any $1<r<p$ there exists a constant $D$ so that for all  
$(a_i)$
and all $M$
$$
{\biggl\Vert \sum_{i=1}^{M} a_i
x_i\biggr\Vert}\le D{\biggl(\sum_{i=1}^M {|a_i|}^r\biggr)}^{1\over  
r}.
$$
\endproclaim

\demo{Proof} We use the convention that ${1\over p}+{1\over q}=1$ and
${1\over r}+{1\over s}=1$.  Given the scalars, assume by
homogeneity that $\sum_{i=1}^M|a_i|^r = 1$. Define
$$
B_n=\{i\bigm|2^{-n-1}<|a_i|\le 2^{-n}\}$$
and write
$$
{\Bigl\Vert \sum_{i=1}^{M}
a_ix_i\Bigr\Vert}
\le
\sum_{n}\Bigl\Vert\sum_{i\in B_n}a_ix_i\Bigr\Vert
\le C \sum_{n} 2^{-n}{|B_n|}^{1\over p}r.
$$
Now just compute:
$$\sum_{n} 2^{-n}{|B_n|}^{1\over p}
=\sum_{n} 2^{-{{nr}\over p}} 2^{n{({r\over p}-1)}}
{|B_n|}^{1\over p}
\le
{\Bigl(\sum_{n} 2^{-nr}{|B_n|}\Bigr)}^{1\over p}
{\Bigl(\sum_{n} 2^{nq{({r\over p}-1)}}\Bigr)}^{1\over q}
$$ 

$$
= 2^{{r\over p}} {\Bigl(\sum_{n} 2^{-(n+1)r}|B_n|\Bigr)}^{1\over  
p}D(p,r)
\le  2^{{r\over p}}{\Bigl(\sum_{i=1}^M |a_i|^r\Bigr)}^{1\over p}  
D(p,r)
$$
$$
= 2^{{r\over p}} D(p,r) {\Bigl(\sum_{i=1}^M |a_i|^r\Bigr)}^{1\over  
r}.
$$

\enddemo
\par
The condition in Theorem 3.3 is rather hard to check.  However, a  
much simpler
criterion is sufficient.
\proclaim{Theorem 3.5} Suppose the dual space of $X$, $X^*$, has type  
$p>1$.
Then for every normalized
weakly null sequence $\{y_n\}$ in $X$ there exists a subsequence and  
a sequence
 $\{x_n\}$ in $X^*$ biorthogonal to it which has an
upper $p$-estimate.  In particular, $X$ satisfies the hypothesis of  
3.3,
and thus is \pnsh\ for all $N>p$.
\endproclaim
\par
In view of the fact that every space with non-trivial type also
has a dual with some non-trivial type, we can state the following  
corollary.

\proclaim{Corollary 3.6}
Suppose $X$ has type.  Then for some $N$, $X$ is \pnsh.\endproclaim
\demo{Proof (of 3.5)}
Recall that the fact that $X^*$ has type $p$ means that there
is a constant $T_p$ so that

$$
\int_0^1{\biggl\Vert \sum_{i=1}^M r_i(t)
x_i\biggr\Vert} dt \le T_p{\biggl(\sum_{i=1}^M
{\Vert x_i\Vert}^p\biggr)}^{1\over p}\qquad\forall M\in {\cal N}
$$
for any finite number of elements $x_1,\dots x_n$ (where $r_i$ are
the Rademacher functions).  Now suppose that
$(y_n)$ is a normalized weakly null sequence, which we can assume is
basic by passing to a subsequence.  We can find a bounded
sequence $(x_n)$ in $X^*$ of functionals biorthogonal to $(y_n)$.  
Since $l_1$ i
s not embeddable in $X^*$ we know by Rosenthal's theorem that we can
find a weakly
Cauchy subsequence of $(x_n)$.  Pass to the odd terms of $(y_n)$,  
relabel and
replace $(x_n)$ with $(x_{2n+1}-x_{2n})$.  Then we have a  
biorthogonal system
 $(x_n,y_n)$ with $(x_n)\to 0$ weakly.  By proposition 3.1 and the  
remark
following it we know that $(x_n)$ has an unconditional spreading  
model $F$.
Now $F$ is finitely representable in $X^*$ (this means that given any  
finite
dimensional subspace of $F$ we can find a $1+\epsilon$--isomorphic  
copy of
that subspace in $X^*$, see again  \cite{B} or  \cite{BS 1}).  Since  
the
definition of type is local, $F$ will also have type $p$ with  
constant
$\le T_p$.  Since the basis $e_n$ of $F$ is unconditional with  
constant
at most 2, it has an upper $p$--estimate with constant less than or  
equal
to $2T_p$. Thus we have satisfied the hypothesis of Theorem 3.3.
\enddemo


\Refs
\widestnumber\key{AAD}


\ref\key AAD
\by {R. Alencar, R. Aron and S. Dineen}
\paper {A Reflexive space of holomorphic functions in infinitely
many variables} \jour Proc. Amer. Math. Soc. \vol 90
\yr 1984 \pages 407-411
\endref


\ref\key ACG
\by R. Aron, B. Cole and T. Gamelin
\paper Spectra of algebras of analytic functions on a Banach space
\toappear
\endref


\ref\key B
\by B. Beauzamy
\paper Banach--Saks properties and spreading models
\jour Math. Scand. \vol 44 \yr 1979 \pages 357--384
\endref




\ref\key BS 1
\by A. Brunel and L. Sucheston
\paper On $B-$convex Banach spaces
\jour Math. Sytems Thy. \vol 7 \yr 1973 \pages 294--299
\endref


\ref\key BS 2
\by A. Brunel and L. Sucheston
\paper On $J-$convexity and some ergodic super-properties of Banach  
spaces
\jour Trans. Amer. Math. Soc. \vol 204 \yr 1975 \pages 79--90
\endref




\ref\key CCG
\by {T. Carne, B. Cole and T. Gamelin.}
\paper {A uniform algebra of analytic functions on a Banach space}
\jour {Proc. Amer. Math. Soc.}
\vol {314}\yr {1989} \pages {639-659}
\endref






\ref\key CGJ
\by B. Cole, T. Gamelin, and W. B. Johnson
\paper Analytic disks in fibers over the unit ball of a Banach space
\toappear
\endref


\ref\key D
\by J. Diestel
\book Sequences and series in Banach spaces
\yr 1984
\publ Springer--Verlag \publaddr New York
\endref










\ref\key F
\by J. Farmer
\paper Polynomial reflexivity in Banach spaces
\toappear
\endref




\ref\key JP
\by J. Jaramillo and  Prieto
\paper Weak-polynomial convergence on a Banach space
\toappear
\endref






\ref\key M
\by J. Mujica
\book   Complex analysis in Banach spaces
\bookinfo Notas de Mathematica, vol. 120
\publ North-Holland \publaddr Amsterdam
\yr 1986
\endref


\ref\key P
\by H. R. Pitt
\paper A note on bilinear forms
\jour J. London Math. Soc. \vol 11  \yr 1936   \pages 174--180
\endref




\ref\key R 1
\by R. Ryan
\book Applications of topological tensor products to infinite
dimensional holomorphy \bookinfo Thesis, Trinity College, Dublin
\yr 1980
\endref




\ref\key R 2
\by R. Ryan
\paper Weakly compact holomorphic mappings on Banach spaces
\jour Pac. J. Math. \vol 131 \yr 1988 \pages 179-190
\endref




\endRefs
\bye